\documentclass[reqno]{amsart}
\usepackage{amsmath}
\usepackage{amsfonts}
\usepackage{epsfig}
\usepackage{mathrsfs}
\usepackage{latexsym}

\newtheorem{theorem}{Theorem}[section]
\newtheorem{lemma}[theorem]{Lemma}
\newtheorem{prop}[theorem]{Proposition}
\newtheorem{cor}[theorem]{Corollary}
\newtheorem*{theorem*}{Theorem}

\theoremstyle{remark}
\newtheorem*{remark*}{Remark}

\numberwithin{equation}{section}

\def\cartanconst{C}
\def\binomconst{B}
\def\hypconst{G}
\def\hypsharpconst{G^\#}
\def\bigeigenconst{\Lambda}
\def\smalleigenconst{\lambda}
\def\bigeigensharpconst{\Lambda^\#}
\def\smalleigensharpconst{\lambda^\#}
\def\binomconst{B}
\def\dufresnoyconst{K}
\def\dufresnoyconstval{12{,}672\left(2.6\cdot10^7\log B + 10^8\right)^{6(4\log B + 20)}\hypconst\hypsharpconst}
\def\covconst{c}

\def\defeq{\stackrel{\mathrm{def}}{=}}
\def\ds{\displaystyle}
\begin{document}
\title[Landau's Theorem for Holomorphic Curves]
{Landau's theorem for holomorphic curves in projective space and
the Kobayashi metric on hyperplane complements}
\author{William Cherry}
\thanks{First author was supported by 
the Humboldt Foundation, NSA grant MDA904-02-1-0093, and MSRI}
\address{Department of Mathematics\\
University of North Texas\\
1155 Union Circle \#311430\\
Denton, TX  76203\\USA}
\email{wcherry@unt.edu}
\author{Alexandre Eremenko}
\thanks{Second author was supported by 
the Humboldt Foundation and NSF grants DMS-0555279 and DMS-0244547}
\address{Department of Mathematics\\
150 N.~University Street\\
Purdue University\\West Lafayette, IN  47907\\
USA}
\email{eremenko@math.purdue.edu}
\date{November 2, 2009}
\begin{abstract} We prove an effective version of a theorem of Dufresnoy:
For any set of $2n+1$ hyperplanes in general position
in $\mathbf{P}^n$, we find an explicit constant $\dufresnoyconst$
such that for every holomorphic map $f$ from the unit disc 
to the complement of these hyperplanes,
we have $f^\#(0)\leq \dufresnoyconst$, where 
$f^\#$ denotes the norm of the derivative
measured with respect to the Fubini-Study metric.

This result gives 
an explicit lower bound on the Royden function,
\textit{i.e.,} the ratio of the Kobayashi metric on the hyperplane
complement to the Fubini-Study metric.
Our estimate
is based on the potential-theoretic method of Eremenko and Sodin.
\end{abstract}
\maketitle
\section{Introduction}
\label{Intro}

Let $H_1, \ldots, H_q$ be hyperplanes in general position
in complex projective space $\mathbf{P}^n,$ $q\geq 2n+1$.
Being in general position
simply means that any $n+1$ hyperplanes have empty intersection.
By a theorem of Dufresnoy
\cite[Th.~VIII]{DufresnoyNormFam},
there exists a constant $K,$ depending only on the 
hyperplanes \hbox{$\{H_1,\ldots,H_q\}$} and the dimension $n,$
such that if $f$ is a holomorphic map from the unit disc
to $\mathbf{P}^n$ which omits the hyperplanes $H_1,\ldots, H_q,$
then $f^\#(0)\le K,$ where $f^\#$ denotes the norm of the
derivative of $f$
with respect to the Fubini-Study metric on $\mathbf{P}^n$
and can be defined by
$$(f^\#)^2= \ds \frac{\ds\sum_{j<k}
        |f_jf_k'-f_kf_j'|^2}{\| f\|^4},
$$
where $\| f\|^2=|f_0|^2+\ldots+|f_n|^2$
for some choice of homogeneous coordinate
functions \hbox{$[f_0:\ldots:f_n]$} for $f.$
In modern terms (see for instance
\cite[\S3.10]{Kobayashi} or \cite[\S VII.2]{Lang}), Dufresnoy's theorem
says that the complement of  $2n+1$ hyperplanes in general position
in $\mathbf{P}^n$ is complete hyperbolic and hyperbolically embedded in
$\mathbf{P}^n.$ 
\par
In \cite[p.~25]{DufresnoyNormFam},
Dufresnoy remarks that the constant $\dufresnoyconst$ depends on
the hyperplanes $H_j$ in a ``completely unknown'' way.
Our purpose is to give an explicit estimate for the constant 
$\dufresnoyconst,$ and therefore an explicit lower bound on the ratio
of the infinitesimal Kobayashi metric (a.~k.~a.\ Royden's function
\cite{Lang}) on the hyperplane complement as
compared to the Fubini-Study metric on projective space.
\par
To that end, let $[X_0:\ldots:X_n]$ be homogeneous coordinates 
on $\mathbf{P}^n,$ and let $H_0,\ldots,H_{n}$ be $n+1$ hyperplanes
in general position in $\mathbf{P}^n$ given by linear defining forms
$$
	H_j(X_0,\ldots,X_n) = a_{j0}X_0+\ldots+a_{jn}X_n
$$
normalized such that
$$
	||H_j||^2 = |a_{j0}|^2+\ldots+|a_{jn}|^2 = 1.
$$
Consider the $(n+1)\times(n+1)$-matrix $A=(a_{jk})$ and let
$$
	0<\lambda_0 \le \lambda_1 \le \ldots \le \lambda_n
$$
be the $n+1$ eigenvalues of the matrix $AA^*$ listed in non-decreasing
order.  We then define the following quantities:
\begin{equation}\label{eigenvaleqn}\begin{split}
	\smalleigenconst(H_0,\ldots,H_n)\defeq \lambda_0 
	&\qquad\textnormal{and}\qquad
	\bigeigenconst(H_0,\ldots,H_n)\defeq \lambda_n\\
	\smalleigensharpconst(H_0,\ldots,H_n)\defeq
	\frac{\sqrt{\lambda_0\lambda_1}}{\lambda_n}
	&\qquad\textnormal{and}\qquad
	\bigeigensharpconst(H_0,\ldots,H_n)\defeq
	\frac{\sqrt{\lambda_{n-1}\lambda_n}}{\lambda_0}.
\end{split}\end{equation}
\par
We can now state our main theorem as
\begin{theorem*}Let $H_0,\ldots,H_{2n}$ be $2n+1$ hyperplanes in
general position in $\mathbf{P}^n.$ Let
$$
	\hypconst=\max_{0\le j_0 < \ldots < j_n\le2n}
	\max \{ \log\bigeigenconst(H_{j_0},\ldots,H_{j_n}),
	\log(n+1)-\log\smalleigenconst(H_{j_0},\ldots,H_{j_n})\},
$$
and
$$
	\hypsharpconst=\min_{0\le j_0 < \ldots < j_n\le2n}
	\frac{1}
	{\smalleigensharpconst(H_{j_0},\ldots,H_{j_n})}.
$$
Let $\binomconst$ be the binomial coefficient
$$
	\binomconst={{2n+1}\choose{n+1}},
$$
and 
$$
	\dufresnoyconst=\dufresnoyconstval.
$$
Let $f$ be a holomorphic map from the unit disc $|z|<1$
to $\mathbf{P}^n$
omitting $H_0,\ldots,H_{2n}.$
Then, 
$$
	f^\#(0)	\le \dufresnoyconst.
$$
\end{theorem*}

When $n=1,$
the classical theorem of Landau \cite{Landau}
gives effective upper bounds, and very good bounds
were obtained in \cite{Hempel} and \cite{Jenkins}.
Effective, though very non-sharp, estimates can also
be derived from Nevanlinna's theory with precise error terms
\cite[\S 5.7--5.8]{CherryYe}.
For some symmetrical arrangements of the omitted points,
one even knows
sharp estimates \cite{BonkCherry}, \cite{BonkCherryTri}.

Dufresnoy's proof of his theorem was based on the deep results
of Bloch and Cartan  \cite{Cartan}, \cite[Ch.~V]{Lang}.
It is not clear whether this approach can give an effective estimate.
One attempt is this direction is \cite{Hall}.

The theorem of Dufresnoy can also be obtained from Borel's theorem \cite{Lang}
by using a compactness argument called
the Zalcman-Brody rescaling lemma
\cite[\S3.6]{Kobayashi}, \cite[\S III.2]{Lang}.
This proof does not give any effective estimate of $\dufresnoyconst.$

\par
If $f$ omits $2^n+1,$ instead of $2n+1,$ hyperplanes,
there is 
an estimate
of Cowen \cite{CowanGeneral} (see also \cite{CowanPlane}),
which he gets by constructing
negatively curved metrics on the hyperplane complement.
It seems difficult to adapt
this method to the case of $2n+1$ omitted hyperplanes, except of course
for $n\le2$ when they are the same.  
For arbitrary $n,$ there is also a negative curvature approach of
Babets \cite{BabetsI},
\cite{BabetsII},  but
Babets's estimates involve higher order derivatives or 
higher associated curves. 

In this paper, we obtain an explicit estimate by putting
in a quantitative form
the
argument from \cite{EremenkoCurves}, where another proof 
of Dufresnoy's theorem was
given. The proof in \cite{EremenkoCurves} was based on potential-theoretic 
considerations which also give a form
of the Second Main Theorem for holomorphic curves
\cite{EremenkoSodin}.
A survey of other results obtained with this method is \cite{ICM}.

To prove our main result we combine the method of \cite{EremenkoCurves}
with two new ingredients. The first one is a quantitative form
of the uniqueness theorem for harmonic functions in the form
of a generalization of the three circles theorem of Hadamard
by Nadirashvili \cite{NadirashviliThreeCirc}.
The second ingredient  is
Theorem~\ref{derivthm}, which gives an estimate
for $f^\#(0)$ in terms of the Fubini--Study area of the image of $f$,
in the case that $f$ omits $n+1$ hyperplanes in general position.
This estimate seems to be new even for $n=1$, 
\cite{EremenkoArea}.
\par
The plan of this paper is as follows.  In section~\ref{linalgsec}
we recall some formulas for the Fubini-Study derivative and some
of its basic properties.  In section~\ref{rieszjensensec}, we
recall the Riesz theorem and Jensen's formula.
In section~\ref{areatoderivativesec}, we show using the Poisson
formula and Harnack's inequality that if a holomorphic map $f$ from
a disc to $\mathbf{P}^n$ omits $n+1$ hyperplanes in general position
and covers finite area $\sigma$
measured in the Fubini-Study metric,
then $f^\#(0)$ can be effectively bounded in terms of $\sigma,$
the dimension $n,$ and the omitted hyperplanes.
In Section~\ref{cartansec}, we recall Cartan's Lemma
in a form we will need.  In Section \ref{nadirashvilisec},
we state Nadirashvili's generalization of the three circles
theorem.  In Section~\ref{coveringsec}, we recall a covering
lemma of Rickman that was also used in \cite{EremenkoCurves}.
Finally, in section~\ref{mainsec}, we put the ingredients together
to give our Landau and Schottky theorems.
\par
\textit{Acknowledgment.} Work on this article began when both
authors visited the Christian Albrechts Universit\"at zu Kiel 
with financial support from the Humboldt Foundation and continued while
the first author was a member at MSRI.  The authors would like to thank
these institutions for their financial support and hospitality.

\section{Fubini-Study Derivatives and Hermitian Linear Algebra}
\label{linalgsec}
If
$$
	H_j(X_0,\ldots,X_n)=a_{j0}X_0+\ldots+a_{jn}X_n
$$
are $n+1$ linear forms defining $n+1$ hyperplanes in general position
in $\mathbf{P}^n$ normalized as in the introduction so that
$$
	||H_j||^2=|a_{j0}|^2+\ldots+|a_{jn}|^2=1,
$$
then for any point $[w_0:\ldots:w_n]$ in $\mathbf{P}^n$, we have
\begin{equation}\label{normeqn}
	\smalleigenconst(H_0,\ldots,H_n) \le
	\frac{\ds\sum_{j=0}^n|H_j(w_0,\ldots,w_n)|^2}
	{\ds\sum_{j=0}^n|w_j|^2}\le
	\bigeigenconst(H_0,\ldots,H_n),
\end{equation}
where $\smalleigenconst$ and $\bigeigenconst$ are defined as in
equation~$(\ref{eigenvaleqn})$.

Let $f=[f_0:\ldots:f_n]$ be a holomorphic map from a domain in 
$\mathbf{C}$ to $\mathbf{P}^n$ given by homogeneous coordinate
functions $f_j$ which are holomorphic without common zeros.
Let $f'=[f_0':\ldots:f_n'].$ We recall the following formulas
for the Fubini-Study derivative $f^\#$ of $f:$
\begin{align}
\label{FSi}\ds (f^\#)^2 &= \ds 
	\frac{\partial^2}{\partial z \partial \bar z}
	\log \sum_{j=0}^{n} f_j\bar f_j \\
\label{FSii}&=
	\ds \frac{\ds\sum_{j=0}^n\sum_{k=0}^n\bar f_j' \bar f_k [f_j'f_k
	-f_j f_k']}{\ds\| f\|^4}
	\;=\;\frac{\ds\sum_{j=0}^n\sum_{k=0}^n|f_j|^2|f_k|^2
	\frac{\bar f_j'}{\bar f_j} \left[\frac{f_j'}{f_j}-\frac{f_k'}{f_k}
	\right]}{\ds\| f\|^4}\\
\label{FSiii} &= \ds \frac{\ds\sum_{j=0}^{n-1}\sum_{k=j+1}^n
	|f_jf_k'-f_kf_j'|^2}{\ds\| f\|^4}\\
\noalign{\vskip 6pt}
\label{FSiv} &= \ds \frac{\ds||f\wedge f'||^2}{\ds||f||^4}.
\end{align}

\begin{prop}\label{changecoordprop}
Let $H_0,\ldots,H_n$ be $n+1$ hyperplanes in general position
in $\mathbf{P}^n$ given by homogeneous forms normalized so that
the coefficient vectors have norm~$1.$
Let $f=[f_0:\ldots:f_n]$ be a holomorphic curve in $\mathbf{P}^n$
with homogeneous coordinate functions $f_0,\ldots,f_n.$
Let $g=[g_0:\ldots:g_n]$ be the holomorphic curve in $\mathbf{P}^n$ given by
$g_j=H_j(f_0,\ldots,f_n).$ Then,
$$
	\smalleigensharpconst(H_0,\ldots,H_n)\le \frac{g^\#}{f^\#}
	\le \bigeigensharpconst(H_0,\ldots,H_n),
$$
where $\smalleigensharpconst$ and $\bigeigensharpconst$ are defined as in
equation~(\ref{eigenvaleqn}).
\end{prop}

\begin{remark*}Only the left-hand inequality involving
$\smalleigensharpconst$ will be used in the sequel.
\end{remark*}

\begin{proof}
With notation as defined preceding equation~(\ref{eigenvaleqn}),
by linear algebra,
$$
	\lambda_0 \le \frac{||g||^2}{||f||^2} \le \lambda_n
	\qquad\textnormal{and}\qquad
	\sqrt{\lambda_0\lambda_1} \le \frac{||g\wedge g'||}{||f\wedge f'||}
	\le \sqrt{\lambda_{n-1}\lambda_n},
$$
and hence the proposition follows from~(\ref{FSiv}).
\end{proof}

\section{The Riesz Theorem and the Jensen Formula}
\label{rieszjensensec}

In this section we recall the Riesz Theorem and the Jensen Formula,
and we derive a corollary of the Jensen Formula that says a 
subharmonic function with large Riesz mass cannot stay close to 
a harmonic function.

\par
Let $v$ be a $C^2$ function on a domain in $\mathbf{C}.$
We recall the differential operator $dd^c$ which can be defined by
$$
	dd^c v = \frac{\partial^2v}{\partial z \partial \bar z}
	\frac{dA}{\pi} = \frac{1}{4}\Delta v \frac{dA}{\pi},
$$
where $dA=dx\wedge dy$ is the Lebesgue area form on $\mathbf{C}.$
One of the advantages of the 
$dd^c$ notation is that
the factors
$2$ and $\pi$ 
do not appear in fundamental formulas.
As usual the operator $dd^c$ is extended to the space $\mathscr{D}'$
of Schwartz distributions. 
When $v$ is subharmonic, $2dd^c v$ defines a (locally finite, positive)
Borel measure known
as the Riesz measure of $v.$

\begin{theorem}[{\textbf{Riesz Theorem}}]\label{rieszthm}
Let $v$ be subharmonic on the unit disc $|z|\le1$ with finite Riesz mass.
Then,
$$
	v(z)\;\;-\int\limits_{|z|<1}\log|z-\zeta|^2dd^cv(\zeta)
$$
is harmonic on $|z|<1.$ 
\end{theorem}

\begin{theorem}[{\textbf{Jensen Formula}}]\label{greenjensen}
Let $v(z)$ be a difference of two subharmonic functions
on $|z|\le R$ and let $0<r<R.$
Then,
$$
	\int_r^R \frac{dt}{t} \int\limits_{|z|\le t}
	dd^c v
	= \frac{1}{2}\left[\int_0^{2\pi}v(Re^{i\theta})
	\frac{d\theta}{2\pi} - \int_0^{2\pi}v(re^{i\theta})
	\frac{d\theta}{2\pi} \right].
$$
\end{theorem}

\begin{cor}\label{jensencor}
Let $v$ be a subharmonic function on $|z|\le R$ and let $u$ be harmonic
on $|z|\le R.$  Let $r<R$ and assume
$$
	(a)\quad\int\limits_{|z|\le r}dd^cv =1 \qquad\textnormal{and}\qquad
	(b)\quad |v-u| < \varepsilon \quad\textnormal{for~}|z|\le R.
$$
Then $\varepsilon \ge \log\ds\frac{R}{r}.$
\end{cor}

\begin{proof}
On the one hand by~(a) and the Jensen Formula,
$$
	\frac{1}{2}\left[\int_0^{2\pi}v(Re^{i\theta})\frac{d\theta}{2\pi}
	-\int_0^{2\pi}v(re^{i\theta})\frac{d\theta}{2\pi}\right]
	=\int_r^R \frac{dt}{t}\int\limits_{|z|\le t}dd^c v\ge 
	\log \frac{R}{r}.
$$
On the other hand by~(b) since $u$ is harmonic,
$$
	\frac{1}{2}\int_0^{2\pi}[v(Re^{i\theta})-v(re^{i\theta})]
	\frac{d\theta}{2\pi}
	\le\frac{1}{2}\int_0^{2\pi}[u(Re^{i\theta})
	-u(re^{i\theta})+2\varepsilon]\frac{d\theta}{2\pi}
	\le \varepsilon. \eqno{\qedhere}
$$
\end{proof}

\section{From an Area Estimate to a Derivative Estimate}
\label{areatoderivativesec}

For a holomorphic curve omitting $2n+1$ hyperplanes in
general position, our eventual goal is to bound $f^\#(0)$ explicitly in
terms of the $2n+1$ given omitted hyperplanes.
What we will first
bound is the integral of $(f^\#)^2$ on a disc centered at the origin.
In this section, we see how to obtain a derivative bound
from such an area bound.
This section is based on \cite{EremenkoArea}
and appears to be new even in dimension one.

\begin{prop}\label{logderivprop}
Let $f$ be a holomorphic function without zeros on $|z|<1$ such
that $|f(z)| < 1.$ Then,
$$
    \left|\frac{f'(0)}{f(0)}\right | \le -\log |f(0)|^2.
$$
\end{prop}

\begin{proof}
Let $r<1.$ Since
$$
    \frac{f'(z)}{f(z)}=\frac{\partial}{\partial z}\log |f(z)|^2,
$$
differentiating the Poisson formula for $\log |f(z)|^2$ gives us
\begin{align*}
    \left | \frac{f'(0)}{f(0)}\right|
    &= \frac{1}{r}\left|\,\int_0^{2\pi}e^{-i\theta}\log|f(re^{i\theta})|^2
    \frac{d\theta}{2\pi}\,\right|\\
    &\le\frac{1}{r}\int_0^{2\pi}-\log|f(re^{i\theta})|^2
    \frac{d\theta}{2\pi}
    =-\frac{1}{r}\log|f(0)|^2
\end{align*}
because \hbox{$|\log|f(z)|^2|=-\log|f(z)|^2.$}  The statement follows
letting $r\to 1.$
\end{proof}

\begin{theorem}\label{derivthm}
Let $f=[f_0:\ldots:f_n]$ be a holomorphic map from the unit disc $|z|<1$
to $\mathbf{P}^n$ which omits the $n+1$ hyperplanes
\hbox{$H_0,\dots,H_n$} in general position.
If
$$
    \int\limits_{|z|<1}(f^\#)^2\,\frac{dA}{\pi} \le \sigma,
$$
then
$$
	f^\#(0) \le \frac{3\sqrt{2}}{\smalleigensharpconst(H_0,\dots,H_n)}
	\left[(2\log2)\sigma + \log(n+1)
	+\log\frac{\bigeigenconst(H_0,\dots,H_n)}{\smalleigenconst(H_0,\dots,H_n)}\right].
$$
\end{theorem}

Before giving a proof, we make some comments.
The above theorem is really only of interest
when $\sigma\ge1.$  When $\sigma<1$ and $n=1,$ a theorem of Dufresnoy
\cite{DufresnoyIsopar} (or see also \cite[Th.~6.1]{HaymanMero} or
\cite[Th.~2.8.3]{CherryYe})
says
$$
    [f^\#(0)]^2 \le \frac{\sigma}{1-\sigma},
$$
even without the assumption that $f$ is zero free. The linear functions
$f(z)=az$ show the Dufresnoy result is sharp in that setting. The same example
shows that one cannot remove the assumption on the omitted hyperplanes
from Theorem~\ref{derivthm}.
Note also that our estimate does not even tend to $0$ as $\sigma\to0.$
\par
We briefly comment on the sharpness of the coefficient
\hbox{$6\sqrt{2}\log2<5.89$} in front of $\sigma.$
Consider the case $n=1$ and a zero free holomorphic function. Let
$$
    f_m(z) = \left(\frac{z-1}{z+1}\right)^m.
$$
Then, $f_m^\#(0)=m,$ and because $f_m$ wraps the disc around the sphere
so as to cover $m$ hemispheres,
$$
    \sigma_m = \int\limits_{|z|<1}(f_m^\#)^2\,\frac{dA}{\pi}
    = \frac{m}{2}.
$$
Hence,
$$
    f_m^\#(0) = 2\sigma_m
$$
no matter how large $m$ is, and
thus the $5.89$ we obtain in front of $\sigma$ is at worst $2.95$
times too big.

\begin{proof}
By first working on a disc of radius $\rho<1$ and then taking a limit
as \hbox{$\rho\to 1,$} we may assume without loss of generality
that $f$ is analytic on \hbox{$|z|\le 1.$}
By abuse of notation, let $H_j$
also denote the linear forms defining the hyperplanes, and normalize
them so that $||H_j||=1.$
Let \hbox{$u_j=\log|H_j\circ f|^2,$}
and without
loss of generality, assume the least harmonic majorant of the $u_j$
on the unit disc $|z|\le 1$ is $0$ and that \hbox{$\max u_j(0)=u_0(0).$}
Let \hbox{$u=\log\sum|f_j|^2$}
and let \hbox{$v(z)=\max_j u_j(z).$}  Then $u$ and $v$ are subharmonic
and satisfy
\begin{align}\label{derivthmstar}
	v&\le \log\sum |H_j\circ f|^2 \le u+\log\bigeigenconst(H_0,\dots,H_n)
	\textnormal{~and}\\
	\label{derivthmstarstar}v&\ge
	\log\sum|H_j\circ f|^2 - \log(n+1)\ge
	u+\log\smalleigenconst(H_0,\dots,H_n)-\log(n+1).
\end{align}
Let $g_j = H_j\circ f$ and $g=(g_0,\dots,g_n).$
Because
$$
	\frac{\ds\sum_{j=0}^n\sum_{k=0}^n|g_j|^2|g_k|^2}
	{\ds\| g\|^4}=1,
$$
we have
$$
	\left(\smalleigensharpconst(H_0,\dots,H_n)f^\#(0)\right)^2 
		\le (g^\#(0))^2
	\le 2\max_j \left|\frac{g_j'(0)}{g_j(0)}\right|^2
	\le 2(u_0(0))^2
$$
by Proposition~\ref{changecoordprop},
formula~(\ref{FSii}), Proposition~\ref{logderivprop} and our assumption that
$$u_0(0)=\max_j u_j(0).
$$
Thus, it suffices to bound $u_0(0).$
\par
Let $0<r<1.$ From the Jensen formula
and the assumed bound on the integral of
$(f^\#)^2,$ we have
$$
    \frac{1}{2}\left[\int_0^{2\pi}u(e^{i\theta})\frac{d\theta}{2\pi}
    -\int_0^{2\pi}u(re^{i\theta})\frac{d\theta}{2\pi}\right]
    =\int_r^1 \frac{dt}{t}\int_{|z|\le t} (f^\#)^2\frac{dA}{\pi}
    \le \sigma \log \frac{1}{r}.
$$
Thus, there is some point $|z_0|$ with $|z_0|=r$ so that
\begin{align*}
	v(z_0)+\log(n+1)-\log\smalleigenconst(H_0,\dots,H_n)\ge u(z_0)&\ge -2\sigma\log \frac{1}{r} +\int_0^{2\pi}
	u(e^{i\theta})\frac{d\theta}{2\pi}\\
	& \ge
	-2\sigma \log \frac{1}{r}-\log\bigeigenconst(H_0,\dots,H_n),
\end{align*}
where the outside inequalities follow from~(\ref{derivthmstar}),
(\ref{derivthmstarstar}), and the fact that $v=0$ when $|z|=1.$
Let $j$ be the index such that
$$
    u_j(z_0)=v(z_0)
$$
and then apply Harnack's Inequality to conclude
\begin{align*}
	u_0(0)\ge u_j(0) &\ge \frac{1+r}{1-r}u_j(z_0)\\
	&\ge
	\frac{1+r}{1-r}\left[-2\sigma\log \frac{1}{r} - \log(n+1)
	-\log\frac{\bigeigenconst(H_0,\dots,H_n)}
	{\smalleigenconst(H_0,\dots,H_n)}\right].
\end{align*}
Taking $r=1/2,$ we get
$$
	u_0(0)\ge-3\left[(2\log 2)\sigma +\log(n+1)+\log
	\frac{\bigeigenconst(H_0,\dots,H_n)}
	{\smalleigenconst(H_0,\dots,H_n)}\right].\qedhere
$$
\end{proof}

\section{Cartan's Lemma}
\label{cartansec}

In this section we recall Cartan's Lemma, 
a well-known fundamental estimate whose
significance in the study of function theory was recognized by Bloch
and made rigorous by Cartan.

\begin{theorem}[{\textbf{Cartan}}]\label{blashthm}
Let $\mu$ be a finite Borel measure on $|z|<1,$ and let
$\Phi$ be the Blaschke potential
$$
	\Phi(z)=\int\limits_{|\zeta|<1}\log\left|\frac{z-\zeta}
	{1-\bar\zeta z\vphantom{T^{T^T}}}
	\right|\,d\mu(\zeta).
$$
Let $0<r<1,$ let $\eta<1/4e,$ and let $|z_0|\le r.$
Then, there is a countable family
of exceptional discs $D_j,$ the sum of whose radii does not exceed
$4e\eta,$ such that for those $z$ with $|z|\le r$ and not in any of the
$D_j,$ we have
$$
	\Phi(z)>\cartanconst(r,\eta)\Phi(z_0),
$$
where
$$
	\cartanconst(r,\eta) = \frac{4}{(1-r)^2}\log\frac{1}{\eta}.
$$
\end{theorem}

Theorem~\ref{blashthm} is stated and proved for finite sums
in Cartan \cite[\S II.21]{Cartan}, which is also
reproduced in \cite[Th.~VIII.3.3]{Lang}.  The argument for 
general $\mu$ with finite mass is the same, so we omit the proof here.

\begin{cor}\label{cartancor}
Let $v$ be a negative subharmonic function on $|z|<1$ such that
for all $\rho<1,$
$$
	\int\limits_{|z|\le\rho}dd^cv<\infty.
$$
Let $r,$ $\eta,$ $z_0$ and $\cartanconst(r,\eta)$ be as in the theorem.
Then,
for all $|z|\le r$ outside of a countable set of discs the sum
of whose radii is at most $4e\eta,$ we have
$$
	v(z) \ge \cartanconst(r,\eta)v(z_0).
$$
\end{cor}

\begin{remark*}Note that the radius of each exceptional disc is at most $r,$
and hence the collective area of the exceptional discs is at most
$4\pi e\eta r.$
\end{remark*}

In \cite[Th.~V]{Cartan},
Cartan states this result in the case
$v=\log|g|$ for an analytic function $g$ with
$|g|<1,$ but with a slightly worse 
expression for $\cartanconst(r,\eta).$ Cartan's proof only makes 
use of Theorem~\ref{blashthm} in the case of a discrete measure,
rather than for a general measure as we shall do.

\begin{proof}
If $v(z_0)=-\infty,$ the estimate is trivial, so we
assume $v(z_0)$ is finite. Let $M$ be a large constant so that
$$
	M\log r <\cartanconst(r,\eta)v(z_0),
$$
and let $\tilde v(z)=\max\{v(z),M\log|z|\}.$ Note that $\tilde v$
is subharmonic,
that $\tilde v(z)=0$ for $|z|=1$ by the assumption that $v(z)<0,$
and that
$$
	\int\limits_{|z|<1}dd^c\tilde v < \infty.
$$
By the Riesz Theorem,
$$
	\tilde v(z)\;-\int\limits_{|z|<1}\log\left|
		\frac{z-\zeta}{1-\bar\zeta z\vphantom{T^{T^T}}}
	\right|^2 \,dd^c\tilde v(\zeta)
$$
is harmonic on $|z|<1.$  However, both $\tilde v$ and the integral potential
vanish identically for $|z|=1,$ and so $\tilde v(z)$ is actually equal to the 
integral.  We then apply the theorem to conclude
the existence of exceptional discs as in the statement of this corollary
such that for those $z$ outside those discs and $|z|\le r,$
$$
	\frac{1}{\cartanconst(r,\eta)}\tilde v(z)\ge 
	\tilde v(z_0) > v(z_0) > M\log r \ge M\log|z|.
$$
We thus conclude that for these same $z,$ we have $\tilde v(z)=v(z).$
\end{proof}

\section{Nadirashvili's Generalization of the Three Circles Theorem.}
\label{nadirashvilisec}

In this section we state a result of Nadirashvili
\cite{NadirashviliThreeCirc} specialized to dimension 2.  
We slightly reformulate Nadirashvili's statement to get explicit
constants convenient for our application.

\begin{theorem}[{\textbf{Three Circles Theorem}}]\label{threecircthm}
Let $u$ be a harmonic function in the unit disc with $|u|<1.$
Let \hbox{$0<r<R<1/2$} and assume
\hbox{$|u(z)|<r^\tau$} for $|z|<r$ and some $\tau>0.$ Then,
\hbox{$|u(z)|<2(2R)^\tau$} for $|z|<R.$
\end{theorem}

\begin{proof}
As in \cite[Lem.~2.1]{KM}, the classical Hadamard three circles
theorem implies,
$$
	\int_0^{2\pi}[u(\rho e^{i\theta})]^2\frac{d\theta}{2\pi}
	\le \rho^{2\tau} \qquad r < \rho < 1.
$$
Thus,
$$
	\int_r^{2R}2\rho\int_0^{2\pi}[u(\rho e^{i\theta})]^2
	\frac{d\theta}{2\pi}d\rho \le 
	\frac{(2R)^{2\tau+2}-r^{2\tau+2}}{\tau+1} \le
	(2R)^{2\tau+2}-r^{2\tau+2}.
$$
Trivially,
$$
	\int_0^{r}2\rho\int_0^{2\pi}[u(\rho e^{i\theta})]^2
	\frac{d\theta}{2\pi}d\rho \le r^{2\tau+2}.
$$
Let $z_0$ be a point with $|z_0|=R$ such that
\hbox{$|u(z)| \le |u(z_0)|$} for \hbox{$|z|\le R.$}
Then, by Jensen's Inequality,
\begin{align*}
	(2R)^{2\tau+2}&\ge\int_0^{2R}2\rho\int_0^{2\pi}[u(\rho e^{i\theta})]^2
	\frac{d\theta}{2\pi}d\rho\\
	&\ge \int_0^R 2t \int_0^{2\pi}[u(z_0+te^{i\theta})]^2
	\frac{d\theta}{2\pi}dt\\
	&\ge \int_0^R2t\left[\int_0^{2\pi}u(z_0+te^{i\theta})
	\frac{d\theta}{2\pi}\right]^2dt = [u(z_0)]^2R^2. \qedhere
\end{align*}
\end{proof}

\begin{prop}[{\textbf{Remez Inequality}}]\label{remez}
Let $D$ be a disc of radius $r$ in $\mathbf{C}$ and let $U$ be
a subset of $D$ with area $\alpha\pi r^2.$  Let $P(x,y)$ be a real polynomial
of degree $n$ on $\mathbf{R}^2$ and assume \hbox{$|P(z)|<\varepsilon$} for
$z$ in $U.$ Then, for all $z$ in $D,$
$$
	|P(z)| < \varepsilon\left(\frac{e}{2\alpha}\right)^n.
$$
\end{prop}

\begin{proof} See \cite[Lem.~2]{NadirashviliThreeCirc}. \end{proof}

If $a$ is a complex number and $r>0,$ let $D(a,r)$ denote
the disc of radius $r$ centered at $a.$

\begin{theorem}[{\textbf{Nadirashvili}}]\label{nadthm}
Let $u$ be a harmonic function in the unit disc such
that $|u|<1.$ Let \hbox{$0<\rho<1/5.$}
Let $U$ be a subset of
$D(0,\rho)$ with area $\alpha\pi\rho^2>0.$ Let
$$
	r=\frac{1}{2}\left(\frac{\alpha}{36e}\right)^2,
$$
and assume $r<\rho/2.$
If $|u(z)|<\varepsilon<1/9$ for all $z$ in $U,$
then there is a point $z_0$ in $D(0,\rho)$ such that
\begin{equation}\label{NadOne}
	|u(z)| < 3\sqrt{\varepsilon}, \qquad\textit{for~all~}
	z \in D(z_0,r).
\end{equation}
Moreover,
\begin{equation}\label{NadTwo}
	\log |u(z)| < \frac{\log(9\varepsilon)\log(1/(5\rho))}
	{2\log\left[\frac{8}{5}\left(\frac{36e}{\alpha}\right)^2\right]}
	+\log 2
	\qquad\textit{for~all~}|z|<\rho.
\end{equation}
\end{theorem}

\begin{proof}
Let $z_0$ be a point in $D(0,\rho)$ such that
$$
	\frac{\mathrm{Area}(U\cap D(z_0,r))}{\pi r^2} \ge
	\frac{\alpha}{72}.
$$
Such a point exists because $D(0,\rho)$ is covered by at most
$$
	\left(3\sqrt{2}\left\lceil\frac{\rho}{r}\right\rceil\right)^2
$$
discs of radius $r$ with centers on a square lattice.
Let $n=\lfloor\log_{2r}\varepsilon\rfloor$ and let $P$ be the $n$-th
Taylor polynomial of $u.$ Because \hbox{$D(z_0,r)\subset D(0,1/2),$}
by estimating the error in Taylor's theorem using $|u|<1,$ we know that
for $z$ in $D(z_0,r),$ we have
\begin{equation}\label{tayloreqn}
	|u(z)-P(z)| \le 2^{n+1}r^{n+1} \le
	(2r)^{\ds\log_{2r}\varepsilon}=\varepsilon.
\end{equation}
Thus, $|P(z)|<2\varepsilon$ for $z$
in $D(z_0,r)\cap U.$  By Proposition~\ref{remez}, if
$z$ is in $D(z_0,r),$ then
\begin{align*}
	|P(z)| &\le 2\varepsilon\left(\frac{e\pi r^2}
	{2\mathrm{Area}(U\cap D(z_0,r))}\right)^n\\
	&\le 2\varepsilon(36e/\alpha)^{
	\ds\log_{2r}\varepsilon}=2\cdot\varepsilon^{\ds 1+\log_{2r}(
	36e/\alpha)}=2\sqrt{\varepsilon},
\end{align*}
where the last equality follows from our choice of $r.$
This together with~(\ref{tayloreqn}) gives~(\ref{NadOne}).

Let $v(z)=u(\frac{4}{5}z+z_0).$ Then, $v$ is harmonic on 
the unit disc, $|v|<1,$ and 
$$
	|v(z)| < 3\sqrt{\varepsilon} \qquad\textnormal{for~}
	|z|<\frac{5}{4}r.
$$
Let $\tau$ be such that $3\sqrt{\varepsilon}=(\frac{5}{4}r)^\tau.$
Then, we conclude by Theorem~\ref{threecircthm} that
$$
	|v(z)| < 2(5\rho)^\tau \qquad\textnormal{for~}
	|z| < \frac{5}{2}\rho.
$$
Because $D(0,\rho)\subset D(z_0,2\rho),$ we conclude that
$$
	|u(z)| < 2(5\rho)^\tau \qquad\textnormal{for~}|z|<\rho,
$$
which gives~(\ref{NadTwo}).
\end{proof}

\section{Rickman's Covering Lemma}
\label{coveringsec}

This section describes a covering lemma of
Rickman \cite{Rickman} as corrected in \cite{EremenkoCurves}.
We provide a proof here to get a slightly better constant.
If $a$ is a complex number and $r>0,$ let $D(a,r)$ denote
the disc of radius $r$ centered at $a.$
For a complex number $b$ with $|b|<1$ and an integer $m\ge1,$
denote by
$$
	\rho_m(b)\defeq\frac{1}{2^{m+1}}(1-|b|).
$$

\begin{prop}\label{covprop}
Let $a$ be a complex number with $|a|<1$ and for $j=0\ldots14,$ let
$$
	b_j=a+\frac{3}{2}\rho_m(a)\frac{a}{|a|}e^{2\pi i j/15},
$$
where $a/|a|$ is understood to be $1$ if $a=0.$ Then
$$
	D(a,2\rho_m(a))\subset D(a,\rho_m(a))\;\cup\;\bigcup_{j=0}^{14} 
	D(b_j,\rho_m(b_j)).
$$
\end{prop}

\begin{figure}[htb]
\epsfig{file=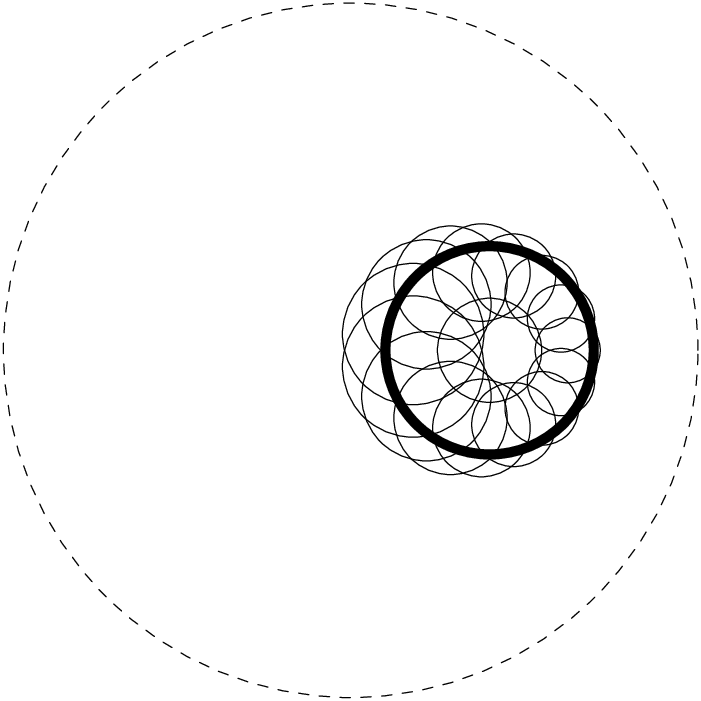,scale=0.7}\hspace{1in}
\raisebox{0.3in}{\epsfig{file=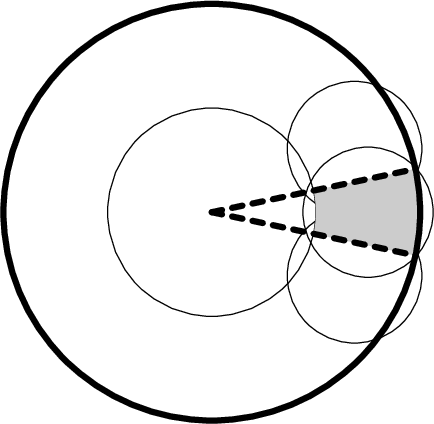,scale=0.8}}
\caption{Illustration of Proposition~\ref{covprop}.}\label{circfig}
\end{figure}

\begin{proof}
Without loss of generality, assume $a\ge0.$
The proof is illustrated by Figure~\ref{circfig}.  The dotted circle
in the left-hand picture represents the boundary of $D(0,1),$  and
the bold circle
is the boundary of $D(a,2\rho_m(a)).$  The picture on the right is
a magnification of the bold circle on the left with only
three of the covering discs of the form $D(b_j,\rho_m(b_j))$ drawn in.
\par
The point is that all the covering discs have radius
at least as big as $\rho_m(b_0).$ Thus, it suffices to show that
the sector of the ring (shaded in the figure on the right)
$$
	\{z : |z-a|>\rho_m(a), |z-a|<2\rho_m(a), 
	-2\pi/30 < \mathrm{arg}(z-a) < 2\pi/30 \}
$$
is contained in the disc $D(b_0,\rho_m(b_0)).$  This is true because
the four corner points of the annular sector are easily seen to be 
contained in 
$D(b_0,\rho(b_0))$  using the Law of Cosines
on the triangles that have as vertices $a,$ $b_0,$
and one of the corner points of the sector, noting that the side lengths
of the triangle and $\rho(b_0)$ vary in direct proportion with
$1-|a|$ as $a$ varies.
\end{proof}

\begin{lemma}[{\textbf{Rickman covering lemma}}]\label{rickmanlemma}
let $\mu$ be a Borel measure on
$|z|<1$ such that $|z|<1$ has finite $\mu$-measure.
Let $m\ge1$ be an integer and let $\covconst>1.$
Then, there is a complex number $a$ with $|a|<1$ such that
$$
	\mu(D(a,2^m\rho_m(a))) \le 16^m\covconst
	\mu(D(a,\rho_m(a)))
$$
and
$$
	\mu(D(0,1/2^{m+1})) \le \covconst\mu(D(a,\rho_m(a))).
$$
\end{lemma}

\begin{proof}
Choose $a$ with $|a|<1$ such that
$$
	\mu(D(a,\rho_m(a))) \ge 
	\frac{1}{\covconst}\sup_{|z|<1}\mu(D(z,\rho_m(z))),
$$
which immediately implies the second inequality of the lemma by considering
$z=0.$ The first inequality follows by iterating Proposition~\ref{covprop}.
\end{proof}

\section{Landau and Schottky Theorems}\label{mainsec}
We now prove our upper bound on the Fubini-Study area covered by
a holomorphic curve omitting $2n+1$ hyperplanes in general position
in $\mathbf{P}^n.$

\begin{theorem}\label{mainthm}
Let $f$ be a holomorphic map from the unit disc $|z|<1$
to $\mathbf{P}^n$
omitting $2n+1$ hyperplanes $H_0,\ldots,H_{2n}$ in general position.
Then, 
$$
	\int\limits_{|z|<\frac{1}{32}}(f^\#)^2\frac{dA}{\pi}
	\le 36\left(2.6\cdot10^7\log B + 10^8
	\right)^{6(4\log B + 20)}\hypconst,
$$
where $\binomconst$ is the binomial coefficient
$$
	\binomconst={{2n+1}\choose{n+1}}
$$
and
$$
	\hypconst=\max_{0\le j_0 < \ldots < j_n\le2n}
	\max \{ \log\bigeigenconst(H_{j_0},\ldots,H_{j_n}),
	\log(n+1)-\log\smalleigenconst(H_{j_0},\ldots,H_{j_n})\},
$$
where $\smalleigenconst$ and $\bigeigenconst$ are defined in 
equation~$(\ref{eigenvaleqn})$.
\end{theorem}

\begin{proof}
By a standard limiting argument we may assume $f$ is holomorphic in
a neighborhood of $|z|\le1.$
Let $[f_0:\ldots:f_n]$ be projective coordinate functions for $f,$
let 
$$
	u=\log \left(|f_0|^2+\ldots+|f_n|^2\right),
$$
and for $j=0,\ldots,2n,$
let $u_j=\log|H_j\circ f|^2,$ where as usual we also use $H_j$
to denote the linear forms defining the hyperplanes normalized so
$||H_j||=1.$  For $\covconst>1,$ by Lemma~\ref{rickmanlemma},
there exists $a$ and $r$
such that
$$
	\int\limits_{|z-a|\le 16r}\!\!\!\!\!\!dd^cu \le
	16^4\covconst\!\!\!\!\!\!\int\limits_{|z-a|\le r}\!\!\!\!\!\!
	dd^cu
$$
and such that
\begin{equation}\label{onequarter}
	\sigma_{a,r}\defeq\!\!\!\!
	\int\limits_{|z-a|\le r}\!\!\!\!
		dd^cu
	\ge \frac{1}{\covconst}
	\int\limits_{|z|\le\frac{1}{32}}dd^cu.
\end{equation}
Let $\tilde u(z)$ be the least harmonic majorant of
$\max\limits_{0\le j\le2n}u_j$ on $|z-a|\le 16r.$  Define
$$
	v(z)=\frac{1}{\sigma_{a,r}}
	[u(a+rz)-\tilde u(a+rz)] \qquad\textnormal{and}\qquad
	v_j(z)=\frac{1}{\sigma_{a,r}}[u_j(a+rz)-\tilde u(a+rz)].
$$
Then $v_j$ are harmonic and negative on $|z|\le 16$
with $\max v_j(z)=0$ when $|z|=16.$ The function $v$ is subharmonic 
on $|z|\le16$ with
\begin{equation}\label{mainstar}
	\int\limits_{|z|\le 1} dd^cv=1
	\qquad\textnormal{and}\qquad
	\int\limits_{|z|\le 16} dd^cv\le16^4\covconst.
\end{equation}
\par
By equation~\ref{normeqn},
for any distinct $n+1$ indices $j_0,\ldots,j_n,$ we have
$$
	\max u_{j_k} \le \log\sum_{k=0}^n|H_{j_k}\circ f|^2
	\le u + \log \bigeigenconst(H_{j_0},\ldots,H_{j_n})
$$
and
$$
	u + \log \smalleigenconst(H_{j_0},\ldots,H_{j_n}) \le
	\log\sum_{k=0}^n|H_{j_k}\circ f|^2 \le
	\max u_{j_k} + \log(n+1).
$$
Thus,
\begin{equation}\label{mainstarstar}
	|v-\max v_{j_k}| \le \varepsilon \qquad
	\textnormal{for~}|z|\le16 \quad\textnormal{where~}
	\varepsilon=\frac{\hypconst}{\sigma_{a,r}}.
\end{equation}
Because the $v_j$ are negative, this implies
$$
	v(z)\le\varepsilon \quad\textnormal{for~}|z|\le16 \qquad
	\textnormal{and}\qquad
	v(z)\ge-\varepsilon \quad\textnormal{for~}|z|=16.
$$
\par
Let $J$ be an index set of cardinality $n+1,$ and let
$$
	U_J=\left\{z : |z|\le2 \textnormal{~and~}
	|v(z)-v_j(z)|<\varepsilon \textnormal{~for~all~}
	j \textnormal{~in~}J\right\}.
$$
By $(\ref{mainstarstar}),$
because we have $2n+1$ hyperplanes total, at every point $z_0$
with $|z_0|\le16,$
there are (at least) $n+1$ distinct indices $j_0,\ldots,j_n$ such
that
$$
	|v(z_0)-v_{j_k}(z_0)| \le \varepsilon \qquad
	\textnormal{for~}k=0,\ldots,n, 
$$
and hence there is at least one such index set $J_0$ such that
$$
	\mathrm{Area}(U_{J_0}) \ge \frac{4\pi}{\binomconst}.
$$
By re-ordering
the indices, without loss of
generality, we will assume $J_0$ contains $0,\ldots,n,$ from now on
we will consider only these indices, and we will denote
$U_{J_0}$ simply by $U.$ By shrinking $U$ if necessary, we may assume
$$
	\mathrm{Area}(U)=\frac{4\pi}{\binomconst}.
$$
\par
By the Jensen Formula and $(\ref{mainstar}),$
$$
	\frac{1}{2}\left[\int_0^{2\pi}v(16e^{i\theta})\frac{d\theta}{2\pi}
	-\int_0^{2\pi}v(2e^{i\theta})\frac{d\theta}{2\pi} \right]
	=\int_2^{16} \frac{dt}{t}\int\limits_{|z|\le t} dd^c v
	\le 16^4\covconst\log 8.
$$
Because $v(z)\ge-\varepsilon$ for $|z|=16,$ we conclude that there is a
point $z'$ with $|z'|=2$ such that 
$$
	v(z')\ge-2\cdot16^4\covconst\log 8 -\varepsilon.
$$
We now apply Corollary~\ref{cartancor} to $v(z)-\varepsilon$ with 
\hbox{$\eta^{-1} > 16e\binomconst$} to conclude that
$$
	v(z)-\varepsilon \ge \cartanconst(1/8,\eta)[v(z')-\varepsilon]
$$
for $|z|\le2$ and outside a collection of discs whose radii $r_j$ are such that
$\sum r_j < 2/\binomconst.$  As $r_j<2,$ this implies
the exceptional discs have collective area less than the area of $U,$
and so there is a point $z_0$ in $U$ with 
$$
	v(z_0)\ge v(z_0)-\varepsilon \ge 
	\cartanconst(1/8,\eta)[v(z')-\varepsilon]
	\ge \cartanconst(1/8,\eta)[-2\cdot16^4\covconst\log8-2\varepsilon].
$$
For $0\le j \le n,$ because $v_j$ is within $\varepsilon$ of $v$ at
$z_0,$ we conclude
$$
	v_j(z_0)\ge-\cartanconst(1/8,\eta)[2\cdot16^4\covconst\log8
							+3\varepsilon].
$$
Applying the Harnack Inequality twice implies
$$
	v_j(z)\ge-9
	\cartanconst(1/8,\eta)(2\cdot16^4c\log8
	+3\varepsilon),\qquad\textnormal{for~}j=0,\ldots,n
	\textnormal{~and~} |z|\le 12.
$$
Letting $\eta^{-1}\to 16e\binomconst,$ we get
for $j=0,\ldots,n$ and \hbox{$|z|\le 12,$}
\begin{equation}\label{maindagger}
	v_j(z)\ge -\frac{2304}{49}
	\log(16e\binomconst)
	(2\cdot16^4\covconst\log8+3\varepsilon).
\end{equation}
\par
Now suppose there is a constant $\delta>0$ such that for any two
indices $j$ and $k$ with $0\le j < k \le n,$ we have
$$
	|v_j(z)-v_k(z)| \le \delta \qquad
	\textnormal{for all~}|z|\le2.
$$
Because at every point $z,$ one of these $n+1$ functions $v_j$
comes within $\varepsilon$ of $v(z)$ and they are all within
$\delta$ of each other, we have that
$$
	|v(z)-v_j(z)|\le \delta+\varepsilon \qquad\textnormal{for~}
	0\le j\le n \textnormal{~and for all~}|z|\le2.
$$
We can now apply Corollary~\ref{jensencor} with $R=2$
and $r=1$ to conclude
$$\varepsilon+\delta\ge\log(2).
$$
We then consider the case that \hbox{$\varepsilon < 1/18,$}
and so from the above,
$$
	\delta>\log 2-\varepsilon>\frac{1}{2}.
$$
\par
In this case, there are indices $j$ and $k$ and a point $\tilde z$
with $|\tilde z|\le 2$ such that
$$
	|v_j(\tilde z)-v_k(\tilde z)|\ge \frac{1}{2}.
$$
Now, $v_j-v_k$ is harmonic on $|z|\le 12,$
$|v_j-v_k|<2\varepsilon$ on $U,$ and
$$
	|v_j-v_k|\le M\defeq 2\cdot\frac{2304}{49}
	\log(16e\binomconst)\left(2\cdot16^4\covconst\log8+\frac{1}{6}\right)
	\quad\textnormal{for~}|z|\le12 \textnormal{~by~}(\ref{maindagger}).
$$
Apply Theorem~\ref{nadthm} to 
$(v_j-v_k)/M,$ which is bounded by $2\varepsilon/M<1/9$
on $U,$ to conclude from~(\ref{NadTwo}) that
$$
	\frac{1}{\varepsilon} \le \frac{36}{M}(2M)^{\ds
	\frac{2\log\frac{8}{5}(36e\binomconst)^2}
	{\log(6/5)}} \le 
	36(2M)^{\ds\frac{2\log\frac{8}{5}(36e\binomconst)^2}
	{\log(6/5)}}.
$$
Letting $\covconst\to1$ and simplifying we get,
$$
	\frac{1}{\varepsilon} < 36\left(2.6\cdot10^7\log B +  10^8
	\right)^{6(4\log B + 20)}.
$$
Because this is much greater than $18,$ there was no harm
in our considering only $\varepsilon < 1/18.$
Combining the above with~(\ref{onequarter}) and recalling
that \hbox{$dd^cu=(f^\#)^2dA/\pi,$} we conclude
$$
	\int\limits_{|z|\le \frac{1}{32}}
	(f^\#)^2\frac{dA}{\pi} \le \sigma_{a,r}
	\le \frac{\hypconst}{\varepsilon},
$$
and hence the theorem follows.
\end{proof}

\begin{theorem}[{\textbf{Landau type theorem}}]\label{landauthm}
Let $H_0,\ldots,H_{2n}$ be $2n+1$ hyperplanes in
general position in $\mathbf{P}^n.$ Let
$$
	\hypconst=\max_{0\le j_0 < \ldots < j_n\le2n}
	\max \{ \log\bigeigenconst(H_{j_0},\ldots,H_{j_n}),
	\log(n+1)-\log\smalleigenconst(H_{j_0},\ldots,H_{j_n})\},
$$
and
$$
	\hypsharpconst=\min_{0\le j_0 < \ldots < j_n\le2n}
	\frac{1}
	{\smalleigensharpconst(H_{j_0},\ldots,H_{j_n})}.
$$
Let $\binomconst$ be the binomial coefficient
$$
	\binomconst={{2n+1}\choose{n+1}},
$$
and let
$$
	\dufresnoyconst=\dufresnoyconstval.
$$
Let $f$ be a holomorphic map from the unit disc $|z|<1$
to $\mathbf{P}^n$
omitting $H_0,\ldots,H_{2n}.$
Then for $|z|<1,$ 
$$
	f^\#(z)	\le \frac{\dufresnoyconst}{1-|z|^2}.
$$
\end{theorem}

When $n=1,$ Landau's original theorem gave a bound on $|f'(0)|$ in terms of 
$|f(0)|$ for analytic functions on the unit disc omitting $0$ and $1.$
In this context,
$$
	f^\#(0) = \frac{|f'(0)|}{1+|f(0)|^2},
$$
and so Theorem~\ref{landauthm} also gives a bound on $|f'(0)|$ in terms
of $|f(0)|.$  Of course, the bound is non-optimal in this special case.
\par
As we discussed in the introduction, the Kobayashi metric on
the complement of the hyperplanes $H_j$ is a metric,
so Theorem~\ref{landauthm}
can be interpreted as a global lower-bound over this hyperplane
complement on the ratio of the
infinitesimal Kobayashi metric to the infinitesimal Fubini-Study metric.
Unfortunately, our method does not give any information in terms of the
geometry of the omitted hyperplanes on where the minimum ratio occurs.
For geometrically symmetric hyperplane configurations, one could hope
to determine these minimum points.  This was done in the one-dimensional
case in \cite{BonkCherry} and \cite{BonkCherryTri}.

Note that the dependence on $\hypconst$ in Theorem~\ref{landauthm}
is essentially best possible.  Indeed, consider $n=1$ and
consider the points $0,$ $m$ and $\infty$ for $m$ large.
Then the map $f_r(z)=e^{rz}$ omits the three points provided
$r<\log m.$  Then, $f_r^\#(0)=r/2$ approaches $(\log m)/2$
as $r$ approaches $\log m.$  Because we have included $0$ and $\infty$
among our three points, $\hypsharpconst=1$ in this case.
A straightforward computation shows that
$$
	\hypconst=\log(2)-\log\left(1-\sqrt{\frac{m^2}{1+m^2}}\right).
$$
An easy application of L'H\^opital's rule also shows that
$$
	\lim_{m\to\infty}\frac{\hypconst}{\log m} = 2,
$$
and hence the linear appearance of $\hypconst$ in the inequality
in Theorem~\ref{landauthm} is correct.  Of course the enormous constant
in front is not optimal.
\par
Now consider the case where the three points are $m,$ $-m,$ and
$\infty$ for $m$ large.  In this case Bonk and
Cherry \cite{BonkCherry} gave a sharp upper bound on $f^\#(0).$
If $f_1$ denotes the universal covering map of
\hbox{$\mathbf{C}\setminus\{1,-1\}$} which sends $0$ to $0,$
then the map $f_m(z)=mf_1(z)$ has the largest spherical derivative at
the origin in comparison with all maps omitting $m,$ $-m$ and $\infty.$
One can even compute using Schwarz triangle functions that
$$
	f_m^\#(0)= m\frac{\Gamma(1/4)^4}{4\pi^2}\approx 4.4m.
$$
One sees as above that here $\hypconst$ is asymptotic to
$2\log m$ as $m\to\infty$ and that  $\hypsharpconst$
is asymptotic to $m,$ and thus, the hyperplane
dependence through the term $\hypsharpconst$ is also not too bad.

The appearance of the combinatorial coefficient
$\binomconst$ in the exponent of our estimate causes the estimate to
deteriorate  rather severely as the dimension $n$ increases.  
It would be interesting to find out
if there must be some dependence on $n$ in the estimate beyond the fact
that $\hypconst$ depends implicitly on $n.$  The appearance of this
combinatorial constant also prevents our method from giving a better
bound when $f$ omits more than $2n+1$ hyperplanes.
If $f$ omits more than $2n+1$ hyperplanes,
the best estimate our method gives is to choose the $2n+1$ hyperplanes that
give the best $\hypconst$ from among the omitted hyperplanes.
This is clearly a deficiency in our approach.  For maps to
$\mathbf{P}^n$ omitting $q$ hyperplanes, it would be
nice to prove an estimate that improves as $q$ increases, as in
\cite[Th.~5.7.4]{CherryYe} when $n=1.$

\begin{proof}
By precomposing $f$ with a M\"obius automorphism of the disc, it
suffices to prove the theorem for $z=0.$

Let $\sigma$ be the upper bound on
$$
	\int\limits_{|z|\le\frac{1}{32}} (f^\#)^2\frac{dA}{\pi}
$$
obtained from Theorem~\ref{mainthm}.
We apply Theorem~\ref{derivthm} taking the best
choice of $n+1$ hyperplanes among our given $2n+1$
omitted hyperplanes to conclude that
$$
    f^\#(0)\le 3\sqrt{2}
	\hypsharpconst[(2\log 2)(32\sigma)+\log(n+1)+2\hypconst]
	\le 11\cdot32\sigma,
$$
since $\log(n+1)+2\hypconst$ is clearly less than $\sigma,$
and the factor $32$ in front of $\sigma$ comes from 
rescaling $D(0,1/32)$ to the unit disc before applying
Theorem~\ref{derivthm}.
\end{proof}

We conclude with a Schottky type theorem which bounds how close
$f(z)$ can get to the omitted hyperplanes depending on the location
of $f(0).$

\begin{theorem}[{\textbf{Schottky type theorem}}]\label{Schottkythm}
With the hypotheses and notation of Theorem~\ref{landauthm}, let
$$
	\delta_j(z)=\frac{|H_j\circ f(z)|}{||f||},
$$
recalling that we have normalized the defining forms of our hyperplanes
so that $||H_j||=1.$ Then for $j=0,\ldots,n$ and $|z|<1,$
$$
	\log\frac{1}{\delta_j(z)}<\frac{1}{1-|z|}\left[
	16\log\frac{1}{\delta_j(0)}+8K^2\right].
$$
\end{theorem}

\begin{proof}
See \cite[\S20]{DufresnoyNormFam}.
\end{proof}

\end{document}